\documentclass[12pt,dvips]{amsart}
\usepackage{amsfonts, amssymb, latexsym, epsfig}
\usepackage{euler, amsfonts, amssymb, latexsym, epsfig, epic,color}

\newtheorem{Theorem}{Theorem}[section] 

\newtheorem{Proposition}[Theorem]{Proposition} 
\newtheorem{Lemma}[Theorem]{Lemma}

\newtheorem{Corollary}[Theorem]{Corollary}

\newtheorem{Definition-Lemma}[Theorem]{Definition-Lemma}
\theoremstyle{remark}

\newcommand\infusion{{\tt infusion}}

\theoremstyle{plain}

\newcommand{\cellsize}{12}
\newlength{\cellsz} 
\newsavebox{\cell}
\sbox{\cell}{\begin{picture}(\cellsize,\cellsize)
\put(0,0){\line(1,0){\cellsize}}
\put(0,0){\line(0,1){\cellsize}}
\put(\cellsize,0){\line(0,1){\cellsize}}
\put(0,\cellsize){\line(1,0){\cellsize}}
\end{picture}}
\newcommand\cellify[1]{\def\thearg{#1}\def\nothing{}%
\ifx\thearg\nothing
\vrule width0pt height\cellsz depth0pt\else
\hbox to 0pt{\usebox{\cell} \hss}\fi%
\vbox to \cellsz{
\vss
\hbox to \cellsz{\hss$#1$\hss}
\vss}}
\newcommand\tableau[1]{\vtop{\let\\\cr
\baselineskip -16000pt \lineskiplimit 16000pt \lineskip 0pt
\ialign{&\cellify{##}\cr#1\crcr}}}
\begin{document}
\pagestyle{plain}

\mbox{}\vspace{-2.0ex}
\title{Cominuscule tableau combinatorics}
\author{Hugh Thomas}
\address{Department of Mathematics and Statistics, University of New Brunswick, Fredericton, New Brunswick, E3B 5A3, Canada }
\email{hugh@math.unb.ca}

\author{Alexander Yong}
\address{Department of Mathematics, University of Illinois at Urbana-Champaign, Urbana, IL 61801}

\email{ayong@uiuc.edu}

\date{October 21, 2014}

\maketitle
\section{Introduction}
The cominuscule Schubert calculus rule of \cite{Thomas.Yong} is based on results of R.~Proctor \cite{Proctor} on poset 
combinatorics, generalizing M.-P.~Sch\"{u}tzenberger's~\cite{Schutzenberger} 
\emph{jeu de taquin} theory.  In this paper, we 
begin with a cominuscule extension of M.~Haiman's \emph{dual equivalence} \cite{Haiman:DE}. One consequence is 
an independent proof of those cases of R.~Proctor's theorem used in \cite{Thomas.Yong}. It also permits us to   
reformulate our rule in a manner that avoids certain arbitrary choices demanded by the original version. 
In addition, we extend S.~Fomin's \emph{growth diagrams} for jeu de taquin to the cominuscule setting and exploit their 
symmetry to give a simple formulation of this case of M.-P.~Sch\"{u}tzenberger's \emph{evacuation involution}. Finally, all of these results and constructions are then used to similarly extend the
$S_3$-symmetric \emph{carton rule} for Littlewood-Richardson coefficients \cite{Thomas.Yong:S3}.  

This work contributes to the theory earlier developed in work of D.~Peterson, R.~Proctor and J.~Stembridge, who show 
that many nice facts for maximal parabolic quotients of the symmetric group hold for 
$d$-complete posets, see, e.g., \cite{Stembridge, Proctor}.  

This paper is entirely combinatorial.  Specifically, we do not discuss the 
geometry that connects this combinatorics to Schubert calculus.
For more on that topic, we refer the reader to \cite{Thomas.Yong} as well as its generalization due to P.~E.~Chaput-N.~Perrin \cite{Chaput.Perrin}. For additional context, 
we also mention that in \cite{Thomas.Yong:K} we extended some of the
combinatorics of this text to the context of $K$-theory. Further research in this direction may be found in, e.g.,
a paper of A.~Buch--V.~Ravikumar \cite{Buch.Ravikumar}, a
joint paper of the authors with E.~Clifford \cite{Clifford.Thomas.Yong},
as well as work of O.~Pechenik \cite{Pechenik}
and of A.~Buch--M.~Samuels \cite{Buch.Samuels}. 

\subsection{Lie-theoretic data and jeu de taquin}
We recall background used in \cite{Thomas.Yong}. This paper centers around
posets associated to seven families of generalized flag manifolds. These
posets are explicitly described on the next page.  
 Although we will present these 
posets in the Schubert calculus 
terminology of our previous work, these posets were 
earlier constructed starting from associated
maximal parabolic subgroups, and called \emph{minuscule posets} in 
\cite{Proctor:euro} (see in particular Section~12 of that paper for
geometric remarks about cohomology of minuscule $G/P$'s). 

Let $G$ be a complex, connected, reductive Lie group
with root system $\Phi$, positive roots $\Phi^{+}$ and
base of simple roots $\Delta$. Fix a choice of maximal
parabolic subgroup $P$ associated to 
a {\bf cominuscule} simple root $\beta(P)$, i.e., if 
$\beta(P)$ occurs in the simple root expansion
of $\gamma\in\Phi^{+}$, it does so with coefficient one. 
Associated to $G$ is the {\bf poset of positive roots} 
$\Omega_{G}=(\Phi^{+},\prec)$ defined by the transitive closure of 
the covering relation $\alpha\prec \gamma$ if $\gamma-\alpha\in \Delta$.
Let
\[\Lambda_{G/P}=\{\alpha\in \Phi^{+}: \mbox{$\alpha$ contains $\beta(P)$ in
its simple root expansion}\}\subseteq \Omega_{G},\] 
the elements of which
we refer to as {\bf boxes}.
Call the lower order ideals of 
$\Lambda_{G/P}$ {\bf straight shapes}, the set of which is denoted by
${\mathbb Y}_{G/P}$.

If $\lambda\subseteq \nu$ are in ${\mathbb Y}_{G/P}$,
their set-theoretic difference is the {\bf skew shape} $\nu/\lambda$.
A {\bf standard filling} of $\nu/\lambda$
is a bijection
\[{\tt label}: \nu/\lambda\to \{1,2,\ldots,|\nu/\lambda|\}
\mbox{ with ${\tt label}(x)<{\tt label}(y)$ whenever $x\prec y$}\] 
(where $|\nu/\lambda|$ denotes the number of boxes of $\nu/\lambda$). This gives a {\bf standard tableau} $T$ of {\bf shape} 
$\nu/\lambda={\tt shape}(T)$. Let ${\rm SYT}_{G/P}(\nu/\lambda)$ be the
set of all such tableaux. 

These tableaux have diagrams similar to those for 
Young tableaux; we now explain this.
The {\bf cominuscule flag varieties} $G/P$ are classified into five
infinite families and two exceptional cases.
For the classical Lie types, we have: 
\begin{itemize}
\item[$A_{n-1}$:] the {\bf Grassmannian} $Gr(k,{\mathbb C}^{n})$,
\item[$B_n$:] the {\bf odd dimension quadric} ${\mathbb Q}^{2n-1}$,
\item[$C_n$:] the {\bf Lagrangian Grassmannian} $LG(n,2n)$,
\item[$D_n$:] the {\bf even dimension quadric} ${\mathbb Q}^{2n-2}$ and\\ the {\bf orthogonal Grassmannian} $OG(n+1,2n+2)$. 
\end{itemize}
The corresponding 
posets $\Lambda_{G/P}$ are the $k\times (n-k)$ rectangle, 
the $1\times (2n-1)$ rectangle, the height $n$ staircase, and a shape with
$2n-2$ boxes, in which all the ranks except the middle one consist of only
a single box. We draw these with
the minimal element in the lower left corner; boxes increase in $\prec$ as we move right or up:
\[\Lambda_{Gr(k,{\mathbb C^n})}:\ \ \ \tableau{{1 }&{ 3 }&{ \ }&{\ }\\{\ }&{ 2 }&{ \ }&{\ }\\
{\ }&{ \ }&{ 4 }&{\ }}, \ \Lambda_{{\mathbb Q}^{2n-1}}:\ \ \ \tableau{{1 }&{2 }&{3 }&{\ }&{\ }}, \ \Lambda_{LG(n,2n)}\cong\Lambda_{OG(n+2,2n+4)}:\ \ \ \tableau{&&&{\ }\\
&&{\ }&{ \ }\\&{3 }&{4 }&{\ }\\
{\ }&{1 }&{2 }&{\ }}\]
\[\Lambda_{{\mathbb Q}^{2n-2}}:\ \ \ \tableau{&&{5 }&{\ }&{\ }& {\ }\\{1 }&{2 }&{3 }&{4 }}\]
We have also inserted standard tableau
of shapes $(3,3,1)/(2,1)$, $(1,1,1)$, $(1,2,2)/(1)$ and $(1,1,2,1)$ 
respectively; we describe a shape as a sequence of column lengths.

For the exceptional Lie types we have: 
\begin{itemize}
\item[$E_6$:] the {\bf Cayley
plane} ${\mathbb O}{\mathbb P}^2$, and 
\item[$E_7$:] the {\bf Freudenthal variety}
$G_{\omega}({\mathbb O}^3, {\mathbb O}^6)$,
\end{itemize}
with posets:
\[\Lambda_{{\mathbb O}{\mathbb P}^2}: \ \ \ \tableau{&&&{\ }&{\ }&{\ }&{\ }&{\ }\\&&&{\ }&{\ }&{\ }\\&&{\ }&{\ }&{\ }\\
{\ }&{\ }&{\ }&{\ }&{\ }}, \ \
\Lambda_{G_{\omega}({\mathbb O}^3, {\mathbb O}^6)}:\ \ \ \tableau{&&&&&&&&{\ }\\&&&&&&&&{\ }\\&&&&&&&&{\ }\\&&&&&&&{\ }&{\ }\\
&&&&{\ }&{\ }&{\ }&{\ }&{\ }\\&&&&{\ }&{\ }&{\ }&{\ }&{\ }\\
&&&&{\ }&{\ }&{\ }\\&&&{\ }&{ \ }&{\ }\\{\ }&{\ }&{ \ }&{\ }&{\ }&{ \ }}\]

Given $T\in {\rm SYT}_{G/P}(\nu/\lambda)$ 
consider $x\in \lambda$,
maximal in $\prec$ subject to the condition that
it is below \emph{some} box of $\nu/\lambda$. 
Associate another standard tableau ${\tt jdt}_{x}(T)$, called the {\bf jeu de
taquin slide} of $T$ into $x$: Let $y$ be the box of $\nu/\lambda$
with the smallest label, among those covering~$x$. Move ${\tt label}(y)$ to $x$, leaving
$y$ vacant. Look for boxes of $\nu/\lambda$ covering $y$ and repeat,
moving into $y$ the smallest label among those boxes covering it. 
Then ${\tt jdt}_{x}(T)$ results when no further 
slides are possible. A {\bf rectification} of $T$
is the result of iterating jeu de taquin slides until 
terminating at a straight shape standard tableau ${\tt rectification}(T)$.
Note that it is not 
obvious that there must be a unique rectification of a tableau; indeed,
this is an important part of the theory which we present here
(Corollary \ref{cor:proctor}).  As an example, one checks that 
there are two possible choices
of orders of slides by which to rectify the $\Lambda_{Gr(k,{\mathbb C}^n)}$
tableau above; using either order, the rectification is $\tableau{{3}\\{1}&{2}&{4}}$.

Given 
$T\in {\rm SYT}_{G/P}(\nu/\lambda)$, consider $x\in\Lambda_{G/P}\setminus
\nu$ \emph{minimal} in $\prec$ subject to being
\emph{above} some element of $\nu/\lambda$. The {\bf reverse jeu de taquin
slide} ${\tt revjdt}_{x}(T)$ of $T$ into $x$ is defined similarly
to a jeu de taquin slide, except we move into $x$ the \emph{largest} of
the labels among boxes in $\nu/\lambda$ covered by~$x$. 

We denote a sequence of slides by the sequence of
boxes $(x_1,\ldots, x_k)$ utilized.

\subsection{Dual equivalence}
We now give a cominuscule extension of M.~Haiman's 
dual equivalence theory \cite{Haiman:DE}: Two tableaux $T$ and $U$ 
are {\bf dual equivalent}, denoted $T\equiv_D U$, if any sequence of slides and reverse slides 
$(x_1,\ldots,x_k)$ for $T$ and $U$ results in tableaux of the same
shape. Clearly, $T\equiv_D U$ implies that ${\tt shape}(T)={\tt shape}(U)$ and moreover, it is easy to prove $\equiv_D$ is an 
equivalence relation on tableaux.

One shape extends another if they can be written as $\nu/\mu$ and $\mu/\lambda$
respectively.  
If $A$ and $B$ are standard tableaux such that ${\tt shape}(B)$
extends ${\tt shape}(A)$, let $A\coprod B$ be
the obvious standard tableau of shape ${\tt shape}(A)\cup {\tt shape}(B)$ 
where the labels of $B$ are increased by $|{\tt shape}(A)|$.

Now suppose that $\lambda \subseteq \mu \subseteq \nu \subseteq \rho$ are
shapes, and let $A$, $B$, $T$, and $U$ be tableaux such that 
\[{\tt shape}(A)=\mu/\lambda, {\tt shape}(T)={\tt shape}(U)=\nu/\mu
\mbox{ and } {\tt shape}(B)=\rho/\nu.\] 
Then it is straightforward
\cite[Lemma~2.1]{Haiman:DE} to show that 
\begin{equation}
\label{eqn:TUpair}
\mbox{if $T\equiv_D U$ then $A\coprod T\coprod B\equiv_D A\coprod U\coprod B$.} 
\end{equation}
Call the replacement of $X:=A\coprod T\coprod B$ by $Y:=A\coprod U\coprod B$ a {\bf Haiman move}. Moreover, call a Haiman move {\bf elementary} if
 the number of boxes $m$ of $T$ and $U$ is:
\begin{table}[h]
\begin{center}
\begin{tabular}{|l||l|l|l|l|l|l|l|}
\hline
$\Phi$ & $A_{n-1}$ & $B_n$ & $C_n$ & $D_n$ & $D_n$ & $E_6$ & $E_7$\\ \hline
$G/P$ & $Gr(k,{\mathbb C}^n)$ & ${\mathbb Q}^{2n-1}$ & $LG(n,2n)$ & ${\mathbb Q}^{2n-2}$ & 
$OG(n+1,2n+2)$ & ${\mathbb O}{\mathbb P}^2$ & $G_{\omega}({\mathbb O}^3, {\mathbb O}^6)$ \\ \hline
$m$ & $3$ & ---   & $4$   & $n$   & $4$   & $5$ & $6$ \\
\hline
\end{tabular}
\end{center}
\end{table}

\noindent
(In ${\mathbb Q}^{2n-1}$, every shape has
exactly one filling, and every dual equivalence class has exactly one
member.)

\begin{Theorem}
\label{thm:main}
For a cominuscule $G/P$: 
\begin{itemize}
\item[(I)] Any two standard fillings of a straight shape $\lambda$ are
dual equivalent. 
\item[(II)] $X, Y\in {\rm SYT}_{G/P}(\nu/\lambda)$
are dual equivalent if and only if they are connected 
by a chain of {\bf elementary Haiman moves}.
\item[(III)] There is a unique straight shape of size $m$ (as given in 
the table above) having  
two standard fillings $T\equiv_D U$. All
other pairs of dual equivalent tableaux of this size are obtained by applying 
a sequence of jeu de taquin slides to this pair. 
\end{itemize}
\end{Theorem}

In \cite{Haiman:DE} the main infinite cases ($Gr(k,{\mathbb C}^n)$, 
$LG(n,2n)$ and $OG(n+1,2n+2)$) of the above theorem were proved. He moreover
notes that many (but not all) 
aspects of his proof generalize to arbitrary posets.
Indeed, our proof of Theorem~\ref{thm:main} follows an approach similar to that
used in his paper. In particular, part (II) is an extension of 
\cite[Proposition~2.4]{Haiman:DE}.

Besides the new exceptional cases of the above theorem, 
our proof differs in two ways from Haiman's. First we 
avoid the need for ``reading word orders'' which were unavailable to us for
the exceptional type cases of our theorem. Second, we introduce a 
simplification (Lemma~\ref{lemma:basic_equiv}) which reduces our proof of 
(III) in the exceptional case 
to a finite check that can be done (tediously) by hand, or,  preferably,
by computer, as is explained in our proof.
This Lemma also simplifies the checks needed
in the previously known cases. We will discuss these aspects in greater detail 
in Section~2.

R.~Proctor~\cite{Proctor} has proved the following corollary in the
greater generality of ``$d$-complete posets'' (not treated here), extending 
\cite{Schutzenberger, sagan, worley}. 
We apply Theorem~\ref{thm:main} to 
obtain an alternative proof for the cominuscule setting.

\begin{Corollary}
\label{cor:proctor}
Given $T\in {\rm SYT}_{G/P}(\nu/\lambda)$, ${\tt rectification}(T)$ is independent
of the order of jeu de taquin slides.
\end{Corollary}

In \cite[p.~5]{Proctor}, R.~Proctor credits D.~Peterson for telling him that
Corollary~\ref{cor:proctor} is true; he writes that Peterson
used a computer to prove this result. In view of our proof
of Theorem~\ref{thm:main}, our proof of this Corollary (for the 
exceptional cases) is also ultimately computationally
based. However, utilizing the technology of dual equivalence allows us to 
avoid checking rectifications of all standard tableaux in types $E_6$ and
$E_7$, and replaces this by a significantly smaller
computer check (small enough to be carried out by hand in type $E_6$).

We now apply dual equivalence to Schubert calculus.
Each $G/P$ is a union of $B_{-}$-orbits whose closures 
$X_w:=\overline{B_{-}wP/P}$ with $wW_P\in W/W_P$
are the {\bf Schubert varieties}. 
The cosets $W/W_P$ correspond bijectively to straight shapes in 
$\Lambda_{G/P}$, 
so we can also refer to the Schubert varieties as $X_\lambda$ for $\lambda\in \mathbb Y_{G/P}$.
(The existence of such a natural correspondence between cosets and order ideals in a poset is a special feature of the cominuscule setting.)  
The
Poincar\'{e} duals $\{\sigma_{\lambda}\}$ of the Schubert varieties form the
{\bf Schubert basis} of the cohomology ring 
$H^{\star}(G/P; {\mathbb Z})$. The {\bf Schubert intersection numbers} $\{c_{\lambda,\mu}^{\nu}(G/P)\}$ 
are defined by
\begin{equation}
\label{eqn:sin}
\sigma_{\lambda}\cdot \sigma_{\mu}=\sum_{\nu\in {\mathbb Y}_{G/P}}c_{\lambda,\mu}^{\nu}(G/P)\sigma_{\nu}.
\end{equation}

If the root system $\Phi$ is not simply-laced, then its roots have two
lengths, referred to as ``long'' and ``short''. If $\Phi$ is simply-laced,
so all roots have the same length, we consider them all to be long.  
Let ${\tt shortroots}(\cdot)$ be 
the number of boxes of a shape that are short roots. 
The following
result relies on our earlier rule \cite[Main Theorem]{Thomas.Yong} to make the
connection to Schubert calculus. 

\begin{Theorem}
\label{thm:nwrl}
For cominuscule $G/P$,
$c_{\lambda,\mu}^{\nu}(G/P)$ equals 
$2^{{\tt shortroots}(\nu/\lambda)-{\tt shortroots}(\mu)}$ times the number of
dual equivalence classes of tableaux of shape $\nu/\lambda$ 
rectifying to a tableau of shape $\mu$. 
\end{Theorem}
Theorem~\ref{thm:nwrl} appears less explicit than our original rule
\cite[Main Theorem]{Thomas.Yong}
(reproduced below as 
Theorem~\ref{thm:originalmain}), although both are 
computationally similar, see the remarks in Section~3.
However, Theorem~\ref{thm:nwrl} has its advantages: it does not 
depend on a fixed choice of tableau of shape $\mu$ to rectify to, 
and its statement is meaningful even in contexts where
Corollary~\ref{cor:proctor} is unavailable.
In this sense, it is more transparent, and possibly useful, e.g., when finding rules for non-(co)minuscule $G/P$. 

In \cite{Thomas.Yong}, a rule was also given for Schubert calculus
for minuscule $G/P$.  Every
minuscule case has a corresponding cominuscule case $(G/P)^{\vee}$ associated
to the Langlands dual group $G^{\vee}$ of $G$.
The Schubert varieties and classes for the minuscule $G/P$
can be indexed by shapes in the corresponding cominuscule $\Lambda_{(G/P)^{\vee}}$. Thus, we also have the following reformulation of the 
minuscule rule of \cite{Thomas.Yong}:

\begin{Corollary} For minuscule $G/P$, $c_{\lambda\mu}^{\nu}(G/P)$ is 
the number of
dual equivalence classes of tableaux in $\Lambda_{(G/P)^\vee}$ 
of shape $\nu/\lambda$ 
which rectify to a tableau of shape $\mu$. 
\end{Corollary}

\subsection{Growth diagrams and their applications}
S.~Fomin's \emph{growth diagrams} provide a way to encode 
jeu de taquin.
In section \ref{ss:comingd}, we explain their straightforward generalization
to cominuscule types.  Growth diagrams make apparent a symmetry of jeu de taquin 
which we refer to 
as the ``infusion involution'' in \cite{Thomas.Yong} (see also 
\cite[Lemma~2.7]{Haiman:DE}).

M.-P.~Sch\"utzenberger defined \emph{evacuation} for 
an arbitrary finite poset.  (See, e.g., the survey \cite{stanley:recent} for background
and references.) Growth diagrams allow us to give a new 
formulation of evacuation for cominuscule posets $\Lambda_{G/P}$.  
As for the classical setting of standard Young tableaux, the fact
that evacuation is an involution is immediate from this perspective.

We refer to shapes of the form $\Lambda_{G/P}/\rho$ as 
{\bf reverse shapes}.  There is a natural bijection between straight 
shapes and reverse shapes,
as follows.  Pick a tableau $T$ of straight shape $\lambda$.
If we apply as many {\tt revjdt} slides as possible to
$T$, the result is a tableau of reverse shape, say $\Lambda_{G/P}/\rho$.    
Since all standard fillings of $\lambda$ are dual equivalent by 
Theorem \ref{thm:main}(I), this shape
only depends on $\lambda$, not on $T$, so we define 
$\lambda^\vee=\rho$.  Since $\Lambda_{G/P}$ is self-dual, the same procedure
can be reversed to go from $\lambda^\vee$ to $\lambda$.  Thus the map
$\lambda \rightarrow \lambda^\vee$ is a bijection. 

For $\lambda,\mu,\nu \in \mathbb Y_{G/P}$ define 
$c_{\lambda,\mu,\nu}(G/P)=c_{\lambda\mu}^{\nu^\vee}(G/P)$.
Because $c_{\lambda,\mu,\nu}(G/P)$ is equal to the number of intersections of 
generic translates by elements of $G$ of the Schubert varieties $X_\lambda$,
$X_\mu$, and $X_\nu$, one
has the obvious $S_3$-symmetries: 
\begin{equation}
\label{eqn:basicsym}
c_{\lambda,\mu,\nu}(G/P)=c_{\mu,\nu,\lambda}(G/P)=c_{\nu,\lambda,\mu}(G/P)
=c_{\mu,\lambda,\nu}(G/P)=c_{\nu,\mu,\lambda}(G/P)=c_{\lambda,\nu,\mu}(G/P).
\end{equation}

In \cite{Thomas.Yong:S3} we 
constructed a {\bf carton rule} for $c_{\lambda,\mu,\nu}$ in the Grassmannian case
that transparently and uniformly explains all of the symmetries~(\ref{eqn:basicsym}).  
As we explain in Section 5, dual equivalence, growth diagrams, and evacuation give us the tools we need
to extend our construction to the cominuscule setting.

\section{Growth diagrams and dual equivalence}
As mentioned above, several steps in our development of cominuscule dual 
equivalence 
will be familiar to readers of \cite{Haiman:DE}. 
However, a crucial step of our argument is different: we 
avoid using ``reading word orders'', which are important 
in \cite{Haiman:DE}, but unavailable to us 
(see further discussion in Section~3).
This necessitates
Lemmas~\ref{lemma:basic_equiv} and~\ref{prop:main}, which 
are deduced in a root-system independent manner.

\subsection{Cominuscule growth diagrams}\label{ss:comingd}

We begin by presenting an extension to the cominuscule setting 
of Fomin's growth diagrams, which 
encode jeu de taquin.  The generalization is
straightforward, but very useful.  
Our proofs parallel 
those in Fomin's Appendix~1 to \cite[Chapter~7]{Stanley}. 

A standard tableau $T$ can be viewed as a {\bf shape chain}, that is to
say, as a sequence of shapes, each successive shape having one more box
than the one before.  
For example, taking $G/P={\mathbb O}{\mathbb P}^2$, we have
\[T=\tableau{&&&{4}\\&&{2}&{3}\\{\ }&{\ }&{\ }&{1}&{5}}\leftrightarrow (1^3) - (1^4) - (1,1,2,1) - (1,1,2,2) - (1,1,2,3) - (1,1,2,3,1),\]
where $(1^3)$ corresponds to the empty boxes of the skew shape, $(1^4)$
gives the shape that also contains the label ``$1$'', $(1,1,2,1)$  is the shape
that contains the labels ``$1$'' and ``$2$'' etc.

One possible rectification sequence of $T$ is given by
\[\tableau{&&&{4}\\&&{2}&{3}\\{\ }&{\ }&{\ }&{1}&{5}} -
\tableau{&&&{\ }\\&&{2}&{4}\\{\ }&{\ }&{1 }&{3}&{5}} -
\tableau{&&&{\ }\\&&{4}&{\ }\\{\ }&{1 }&{2 }&{3}&{5}}- 
\tableau{&&&{\ }\\&&{4}&{\ }\\{1 }&{2 }&{3 }&{5}&{\ }},\]
and each of these skew tableaux has its own shape chain. 
Putting the
shape chain for $T$ atop the shape chains for each of the tableaux
produced in the course of rectifying $T$, we obtain a two-dimensional
array of shapes,
a cominuscule analogue of Fomin's {\bf growth diagram}, which  
in the example at hand is given in Table \ref{table:gd}. 

\begin{table}[h]
\begin{center}
\begin{tabular}{|l|l|l|l|l|l|}
\hline
$(1^3)$ & $(1^4)$ & $(1,1,2,1)$ & $(1,1,2,2)$ & $(1,1,2,3)$ & $(1,1,2,3,1)$\\ \hline
$(1^2)$ & $(1^3)$ & $(1,1,2)$ & $(1,1,2,1)$ & $(1,1,2,2)$ & $(1,1,2,2,1)$\\ \hline
$(1)$ & $(1^2)$ & $(1^3)$ & $(1,1,2)$ & $(1,1,2,1)$ & $(1,1,2,1,1)$ \\ \hline
$\emptyset$ & $(1)$ & $(1^2)$ & $(1^3)$ & $(1,1,2)$ & $(1,1,2,1)$ \\ \hline
\end{tabular}
\end{center}
\caption{A cominuscule Fomin growth diagram\label{table:gd}}
\end{table}

Note that the top row encodes the original tableau, while 
the
left column $\emptyset - (1) - (1^2) - (1^3)$ corresponds
to the tableau $R=\tableau{{1}&{2}&{3}}$, which describes the order of the
jeu de taquin slides in the rectification.  


Growth diagrams can be characterized in the following way:

\begin{Theorem}\label{thm:gd} A rectangular array of straight shapes in $\mathbb Y_{G/P}$
is a growth diagram if and only if for any $2\times 2$ subgrid
\begin{tabular}{|c|c|}\hline
$\alpha$ &  $\beta$\\ \hline
$\gamma$ & $\delta$\\\hline
\end{tabular} the {\bf Fomin growth conditions} hold:
\begin{itemize}
\item[(F0)] $\alpha/\gamma$, $\delta/\gamma$, $\beta/\alpha$, and 
$\beta/\delta$ each consist of a single box;
\item[(F1)] if $\alpha$ is the unique shape containing $\gamma$ and
contained in $\beta$, then $\delta=\alpha$;
\item[(F2)] otherwise there is a unique such shape other than $\alpha$, 
and this shape is $\delta$. 
\end{itemize}
\end{Theorem}
\begin{proof}
We first check that a growth diagram satisfies the growth conditions.
The $i$-th row of the growth diagram defines a tableau $T_i$.  
To verify (F0), consider the $2 \times 2$ subgrid located in rows $i$ and $i+1$, and
columns $j$ and $j+1$.  
Then $\beta/\alpha$ is the position of $j$ in
$T_{i+1}$, while $\delta/\gamma$ is the position of $j$ in $T_i$.
In the course of the jdt slide which changes $T_{i+1}$ to $T_i$, we have that
$\alpha/\gamma$ is the position of the empty box after  
the boxes numbered 1 to $j-1$ have moved 
from their positions in $T_{i+1}$ to their
positions in $T_{i}$, and $\beta/\delta$ is the position of the empty box 
after box $j$ has also moved.  This establishes (F0).  Next, we observe
that condition (F1) is
automatically satisfied given (F0).  (It is included in the growth 
conditions for clarity.)  Finally, if there is a unique shape between
$\beta$ and $\gamma$ other than $\alpha$, then after we have moved
$1$ through $j-1$ from their positions in $T_{i+1}$ to their positions in
$T_i$, then the empty box and the box containing $j$ are not adjacent, with
the result that $j$ occupies the same position in $T_i$ as in 
$T_{i+1}$, which implies (F2).  This establishes that the growth 
conditions hold for each $2\times 2$ subgrid.

Conversely, suppose that we have a rectangular array of shapes satisfying
the growth conditions.  Interpret the leftmost column as a straight shape 
$A$, and interpret the top row as a skew shape $B$.  Now consider the growth diagram 
for the rectification of $B$ in the rectification order given by $A$.  This
new diagram has the same left column and top row as our original diagram,
and both satisfy the growth conditions.  Since the growth conditions suffice
to determine the whole array given the left column and top row, the two arrays
must coincide, and the given array must be a growth diagram.  
\end{proof}

Observe that the Fomin growth conditions are symmetric under a 
transposition about the bottom-left/top-right diagonal.  This leads to an
important tableau-theoretic involution, which we refer to as 
``infusion'' (this is a much older concept, see \cite{Haiman:DE} as well
as, e.g., \cite{Benkart}).
Let $A$ be a standard tableau of shape $\lambda$, and $B$ a standard tableau
of shape $\nu/\lambda$.  We define $\infusion(A,B)=(C,D)$ where $C$ is the
result of rectifying $B$ according to the order given by $A$, and $D$
is the tableau which records the order in which boxes of $\nu$ were emptied
in the rectification procedure.  If we consider the growth diagram for the
rectification of $B$ in the order $A$, the bottom row 
gives the shape chain for $C$, and the rightmost column gives the
shape chain for $D$.  The fact that growth diagrams are transpose 
symmetric then implies that $\infusion(C,D)=(A,B)$; that is to say, that
infusion is an involution.  For future use, we record the notation
that if $\infusion(A,B)=(C,D)$, then $\infusion_1(A,B)=C$ and 
$\infusion_2(A,B)=D$.

In fact, the same proof holds in a slightly more general setting.  
The following fact
was also 
proved in~\cite[Theorem~4.4]{Thomas.Yong}.
\begin{Lemma}
\label{lemma:infusionandrev}
For any standard tableaux $T$ and $U$ such that ${\tt shape}(U)$ extends
${\tt shape}(T)$ then
${\tt infusion}({\tt infusion}(T,U))=(T,U)$.
\end{Lemma}

\subsection{Proof of Theorem~\ref{thm:main}}

Consider the {\bf basic shapes}, 
which are the minimal shapes in each
$\Lambda_{G/P}$ having two standard fillings, as displayed in
Table \ref{basic}.
\begin{table}[h]
\begin{center}
\begin{tabular}{|l||l|l|l|l|l|l|l|}
\hline
$G/P$\!\! & \!$Gr(k,{\mathbb C}^n)$ \!\! & \!\!\!\! ${\mathbb Q}^{2n-1}$\!\! & $LG(n,2n)$ & ${\mathbb Q}^{2n-2}$ \!\!\!\!\! & 
$\! OG(n\!+\! 1,2n\!+\! 2)$\!\! & \!${\mathbb O}{\mathbb P}^2$ \!\!\! & $G_{\omega}({\mathbb O}^3, {\mathbb O}^6)\! $ \\ \hline
$\lambda$ & $(2,1)$ & ---   & $(1,2,1)$   & $(1^{n-3},2,1)$   & $(1,2,1)$   & $(1,1,2,1)$ & $(1,1,1,2,1)$ \\ \hline
& $\tableau{ {\ }\\ {\ } & {\ } \\ \\ }$ & & $\tableau{  &{\ }\\
{\ }&{\ }&{\ }}$ & 
$\tableau{  \\ {\ }&{\ } \\ \\}\begin{array}{c}  \\ \ldots 
\end{array}\tableau {  & {\ } \\ {\ }&{\ }&{\ }\\
\\}$
& 
$\tableau{ &{\ }\\
{\ }&{\ }&{\ }}$
 & $\tableau{  &&{\ } & \\ {\ }&{\ }&{ \ }&{\ }\\ \\}$ &  
$\tableau{  &&&{\ } & \\ {\ } & {\ }&{\ }&{ \ }&{\ }\\ \\}$
\\
\hline
\end{tabular}
\caption{The basic shapes for cominuscule posets.\label{basic}}
\end{center}
\end{table}

We now establish a special case of Theorem~\ref{thm:main}(I), which turns
out to be fundamental:

\begin{Lemma}
\label{lemma:basic_equiv}
For each cominuscule $G/P$, the two tableaux of the
basic shape are dual equivalent. 
\end{Lemma}
\begin{proof}
Consider a sequence of slides 
\begin{equation}
\label{eqn:slideseq}
{\tt revjdt}_{x_1}(\cdot ),\cdots, {\tt jdt}_{x_i}(\cdot ),\cdots
{\tt revjdt}_{x_j}(\cdot ), \cdots
\end{equation}
associated to boxes $\{x_i\}\subseteq \Lambda_{G/P}$.
Let $T_1$ and $T_2$ be the two standard tableaux of the basic shape. 
Call a {\bf direction change} in 
(\ref{eqn:slideseq}) a ${\tt revjdt}$ slide followed by a ${\tt jdt}$ slide, 
or a ${\tt jdt}$ slide followed by a ${\tt revjdt}$ slide. 
We induct on the number of direction changes to show that the 
shapes of $T_1$ and $T_2$ under (\ref{eqn:slideseq}) are the same.

In the base case, there are none, and the conclusion
is a straightforward (but tedious) verification; in the classical types
analyzed in \cite{Haiman:DE}, a similar approach was also suggested. 
However, here our task is actually simpler since we only need to check
for size $m=3$ (in the $Gr(k,{\mathbb C}^n)$ case) and $m=4$ (in the
$LG(n,2n)$ and $OG(n+1,2n+2)$ cases) that reverse slides preserve the
equality of 
shapes of the two tableaux of these sizes. Although this is an infinite check,
the possibilities for how the relative positions of the $m$ boxes can appear
(in relation to $\Lambda_{G/P}$) is small and can be indeed analyzed 
(although we omit the details).   

The other classical types are easy to check.

Finally, in the exceptional types the check is finite.
In type $E_6$, we carried this check out by hand. 
In type $E_7$, the number of cases is
significantly larger.  Though it is still within reach of 
human verification, we preferred to handle this case using 
a simple {\tt Maple} program\footnote{Software available at the authors' websites.} (which we also used to reconfirm our hand-calculations for type $E_6$).  
Our
program constructs all possible (partial) reverse rectifications 
recursively (at each recursive step, each possible reverse jeu de taquin
move is determined). This check takes no more than a few minutes on a 
 computer. 

This concludes the discussion of the base case of this proof.
 
Now we assume that there is at least one direction change.  The first direction
change is of the form ``${\tt revjdt}_{x_c}(\cdot), {\tt jdt}_{x_{c+1}}(\cdot)$''. 
Up until $x_{c}$, we have been solely 
applying ${\tt revjdt}$ slides, obtaining, by the base case, 
$T_1'$ and $T_2'$ of the same
shape. 

Recall that $\Lambda_{G/P}$ is self-dual, and that we refer to a shape of the
form $\Lambda_{G/P}/\lambda$ as a {reverse shape}.

By the base case, there is sequence of slides
\[{\tt revjdt}_{z_1}(\cdot),\ldots, {\tt revjdt}_{z_M}(\cdot)\]
that ``reverse rectify'' $T_1'$ and $T_2'$ to tableaux 
$T_1''$ and $T_2''$ of the same reverse shape. Suppose 
\[{\tt jdt}_{y_M}(\cdot),\ldots, {\tt jdt}_{y_1}(\cdot)\] 
are the 
slides that undo the $\{z_i\}$ slide sequence, 
returning us to $T_1', T_2'$. Observe that
by the self-duality of $\Lambda_{G/P}$ we can interpret 
\[{\tt jdt}_{y_M}(\cdot),\ldots, {\tt jdt}_{y_1}(\cdot), 
{\tt jdt}_{x_{c+1}}(\cdot)\] 
as a sequence of ${\tt revjdt}$ slides for the dual poset to
$\Lambda_{G/P}$ which is, of course, isomorphic to $\Lambda_{G/P}$. 
Finally, concatenating the slides into 
$x_{c+2},\ldots,x_N$ of (\ref{eqn:slideseq}) reintepreted by
${\tt jdt}\leftrightarrow {\tt revjdt}$, we obtain a new sequence
of slides with one fewer direction change that passes through $T_1',T_2'$. 
Thus by induction,
the resulting tableaux have the same shape, and therefore the same
would be true of applying (\ref{eqn:slideseq}) to $T_1,T_2$.
\end{proof}

\begin{Lemma}
\label{prop:main}
Let $c$ and $d$ be two distinct corners of $\lambda\in {\mathbb Y}_{G/P}$. 
There exists a sequence of jeu de taquin slides that
when applied to one of the two standard tableaux of the basic shape $\beta$,
the entry $|\beta|-1$ is sent to $c$ and the entry $|\beta|$ to $d$, while for 
the other standard tableau, $|\beta|$ is sent to $c$ and 
$|\beta|-1$ to $d$.
\end{Lemma}
\begin{proof}
Mark the two corners (i.e., those containing
$|\beta|-1$ and $|\beta|$) of $\beta$ with a ``$\star$''.  We wish to show 
that there is a sequence of jeu de taquin slides moving the two $\star$'s to
$c$, $d$, without ever producing a situation where the two $\star$'s are trying
to move into the same box. Clearly, 
such a sequence of jeu de taquin slides can be constructed by moving each of the 
$\star$'s along the boundary (as drawn in the plane) of the minimal
straight shape containing $c$ and $d$; one
takes the northwest boundary and the other the southeast boundary.   
\end{proof}

\begin{Lemma}
\label{lemma:new} Let $c,d$ be two distinct corners of 
$\lambda \in {\mathbb Y}_{G/P}$.  Then there exist two tableaux $S_1$ and 
$S_2$ of shape $\lambda$, related by a single elementary Haiman move,
such that $S_1$ has $|\lambda|$ in $c$ while $S_2$
has $|\lambda|$ in $d$. 
\end{Lemma}

\begin{proof}
Start with the two fillings of the basic shape. 
By Lemma~\ref{lemma:basic_equiv} these are 
dual equivalent. 
Apply the sequence of slides constructed
in Lemma~\ref{prop:main}.  The result is two dual equivalent 
fillings $B_1$, $B_2$ of a shape 
$\lambda/\gamma$ for some $\gamma$, one having its maximum entry in
$c$, the other having its maximum entry in $d$.  Let $A$ be an arbitrary
standard filling of $\gamma$.  Then $S_i:=A\amalg B_i$ satisfy the 
statement of the lemma.
\end{proof}

For use below, we point out the following facts which follow immediately from
the definition of dual equivalence:
\begin{Lemma}
\label{lemma:clearly}
If $T\equiv_{D} U$ then ${\tt jdt}_x(T)\equiv_D {\tt jdt}_{x}(U)$
and ${\tt revjdt}_x(T)\equiv_D {\tt revjdt}_{x}(U)$.

If $T\equiv_D U$ by an elementary Haiman move, then the same is true
for the tableaux resulting from applying the same slide to $T$ and $U$.
\end{Lemma}

\noindent{\it Conclusion of the proof of Theorem~\ref{thm:main}:}
To prove (I), we induct on $|\lambda|$. The base
case $\lambda=\emptyset$ is obvious. Now suppose $\lambda$ has 
at least one box. By induction, for any corner $c$ of $\lambda$, 
there are (elementary) Haiman moves connecting those tableaux having
$|\lambda|$ in $c$.  Thus we are done if there is only one corner
of $\lambda$, so suppose there are at least two corners $c$ and $d$, and
$T_1, T_2\in {\rm SYT}_{G/P}(\lambda)$ where $T_1$ has $|\lambda|$ in $c$
and $T_2$ has $|\lambda|$ in $d$. By Lemma~\ref{lemma:new}, 
there is 
\[S_1\equiv_D S_2 \mbox{ with } 
{\tt shape}(S_1)={\tt shape}(S_2)=\lambda\] 
such that $S_1$ has $|\lambda|$ in $c$ and
$S_2$ has $|\lambda|$ in $d$. Thus 
\[T_1\equiv_D S_1 \equiv_D S_2 \equiv_D T_2\] 
as desired.

For (II), ``$\Leftarrow$'' is trivial. Conversely, let 
$T\in {\rm SYT}_{G/P}(\lambda)$. By (I), there is a chain of 
elementary Haiman moves 
\[{\tt infusion}_1(T,X)=C_0\equiv_D C_1\equiv_D \cdots
\equiv_D C_N={\tt infusion}_1(T,Y).\] 
Since $X\equiv_D Y$,
${\tt infusion}_2(T,X)={\tt infusion}_2(T,Y)$. Let 
\[D_i={\tt infusion}_2(C_i, {\tt infusion}_2(T,X)).\] 
Then by Lemma~\ref{lemma:clearly}, it follows
that 
\[X=D_0\equiv_D D_1\equiv_D \cdots\equiv_D D_N=Y\] 
is a chain of elementary Haiman moves. 

For (III), the assertion that there are only two fillings of size $m$ is
obvious. That these two fillings are dual equivalent is Lemma~\ref{lemma:basic_equiv}. The second claim follows
by choosing a rectification sequence for the given dual equivalent tableaux.
Since the two resulting tableaux must be different fillings of the
same straight shape, the result follows by the first assertion.\qed

\subsection{Proof of Corollary~\ref{cor:proctor}}
Let $T$ be a skew tableau 
in ${\rm SYT}_{G/P}(\nu/\lambda)$, and write $x_i$ for the box of $T$ with
entry $i$. Let  $A, B \in {\rm SYT}_{G/P}(\lambda)$ encode two possible
rectification orders for $T$.  
Since $A\equiv_D B$, we have that we have that 
$\infusion_1(A,T)=\infusion_1(B,T)$.  
%
%
%
%
\qed

\subsection{Proof of Theorem~\ref{thm:nwrl}}

A pair of tableaux $T,U$ are {\bf jeu de taquin equivalent} if
\[{\tt rectification}(T)={\tt rectification}(U).\] 
They are merely {\bf shape equivalent} if  
\[{\tt shape}({\tt rectification}(T))={\tt shape}({\tt rectification}(U)).\]

\begin{Proposition}
\label{prop:rowcol}
Fix a shape $\nu/\lambda\subseteq \Lambda_{G/P}$. Within each shape
equivalence class, each jeu de taquin equivalence class meets each dual
equivalence class in a unique $T\in {\rm SYT}_{G/P}(\nu/\lambda)$.
\end{Proposition}
\begin{proof}
Fix a choice of $U\in {\rm SYT}_{G/P}(\lambda)$.  
We must show that for any 
$A,B\in {\rm SYT}_{G/P}(\nu/\lambda)$ that are shape equivalent,
there exists a unique $T\in {\rm SYT}_{G/P}(\nu/\lambda)$ such that
${\tt infusion}_{1}(U,A)={\tt infusion}_{1}(U,T)$ (i.e., $T$ and $A$
are in the same jeu de taquin class) and $T\equiv_{D} B$. 

Notice that in fact, if we write $R$ for $\infusion_2(U,B)$, then 
\[{\tt infusion}_{1}(U,\cdot) \mbox{ and }{\tt infusion}_{2}(\cdot,R)\]
are mutually inverse bijections between 
\[\mbox{the dual equivalence class of $B$
and ${\rm SYT}_{G/P}({\tt shape}({\tt infusion}_{1}(U,B)))$.}\] 
Therefore
\[T={\tt infusion}_{2}({\tt infusion}_{1}(U,A),R)\] 
does the job.
\end{proof}

Theorem~\ref{thm:nwrl} then follows immediately from Proposition~\ref{prop:rowcol} and

\begin{Theorem}(\cite[Main Theorem]{Thomas.Yong})
\label{thm:originalmain}
For cominuscule $G/P$, 
let $\lambda,\mu,\nu\in {\mathbb Y}_{G/P}$ and fix 
$T_{\mu}\in {\rm SYT}_{G/P}(\mu)$. 
Then $c_{\lambda,\mu}^\nu(G/P)$ is 
$2^{{\tt shortroots}(\nu/\lambda)-{\tt shortroots}(\mu)}$ times the 
number of standard tableaux of shape $\nu/\lambda$ whose rectification
is $T_\mu$.  

For minuscule $G/P$, let 
$\lambda,\mu,\nu\in {\mathbb Y}_{(G/P)^\vee}$ and fix 
$T_{\mu}\in {\rm SYT}_{G/P}(\mu)$.
Then $c_{\lambda,\mu}^\nu(G/P)$ is 
the 
number of standard tableaux of shape $\nu/\lambda$ whose rectification
is $T_\mu$. \qed
\end{Theorem}

\section{Further discussion of dual equivalence}

\subsection{Computing $c_{\lambda,\mu}^{\nu}(G/P)$} 
Consider $G/P={\mathbb O}{\mathbb P}^2$, the Cayley
plane associated to the root system $E_6$, 
and the skew shape $\nu/\lambda=(1,1,2,3,1)/(1,1,1)$. The seven fillings
are given in Table~\ref{table:sjd}. In this {\bf Haiman table}, the 
rows give the jeu de taquin
equivalence classes, and the columns give the dual equivalence classes,
in agreement with Proposition~\ref{prop:rowcol}.
The rightmost column computes the common {\tt rectification} of the
tableaux in a given row. 
Theorem~\ref{thm:nwrl} says, e.g., that $c_{(1,1,1),(1,1,2,1)}^{(1,1,2,3,1)}({\mathbb O}{\mathbb P}^2)=3$ by counting the middle three columns. Meanwhile
Theorem~\ref{thm:originalmain} says count the three tableaux 
in \emph{either} the second or third row.

In practice, both rules are similar: 
in using Theorem~\ref{thm:originalmain},
we do not know of any general way to avoid essentially checking all
skew tableaux of shape $\nu/\lambda$. So, we basically 
produce much of the information needed to construct a Haiman table,
which encodes all coefficients 
$c_{\lambda,\gamma}^{\nu}(G/P)$ as $\gamma$ varies. (To determine
if two tableaux are dual equivalent, check if one tableau's
rectification sequence works for the other, and produces the same shape.)

\subsection{The Haiman table and the generalized
Robinson-Schensted correspondence}
Organizing one's
thoughts about Schubert intersection numbers this way can be illuminating.
For example, when 
\[G/P=Gr(k,\mathbb{C}^n) \mbox{ and } 
\nu/\lambda=(k,k-1,\ldots,3,2)/(k-1,k-2,\ldots,2,1)\] 
the standard fillings 
are in obvious bijection with the symmetric group $S_k$. The last column is the
``insertion tableau'' of the Schensted insertion algorithm. The 
``recording tableau'' of his algorithm labels the columns. Viewed this way, Proposition~\ref{prop:rowcol} generalizes Robinson-Schensted to arbitrary standard (cominuscule) tableaux, extending an observation of~\cite{Haiman:DE}.
 
\begin{table}[t]
\begin{center}
\begin{tabular}{|l|l|l|l||l|}
\hline
$\tableau{\\&&&{4}\\&&{1}&{3}\\{\ } & {\ }& {\ }&{2}&{5}\\ \\}$ & & & & 
$\tableau{\\ &&&{\ }\\&&{\ }&{ \ } \\{1}&{2}&{3}&{4}&{5}\\ \\}$ \\ \hline
& $\tableau{\\&&&{4}\\&&{2}&{3}\\{\ }& {\ } &{\ } &{1}&{5}\\ \\}$ & 
$\tableau{\\&&&{5}\\&&{1}&{4}\\{\ }&{\ }&{\ }&{2}&{3}\\ \\}$ &
$\tableau{\\&&&{5}\\&&{2}&{4}\\{\ }&{\ }&{ \ }&{1}&{3}\\ \\}$ & 
\tableau{\\ &&&{\ }\\ &&{4}&{\ }\\{1}&{2}&{3}&{5}&{\ }\\ \\}\\ \hline
& $\tableau{\\&&&{5}\\&&{2}&{3}\\{\ } & {\ } &{ \ }&{1}&{4}\\ \\}$ & 
$\tableau{\\&&&{5}\\&&{1}&{3}\\{\ }&{\ }&{ \ }&{2}&{4}\\ \\}$ &
$\tableau{\\&&&{5}\\&&{3}&{4}\\{\ }&{\ }&{\ }&{1}&{2}\\ \\}$ & 
\tableau{\\ &&&{\ } \\ &&{5}&{\ }\\{1}&{2}&{3}&{4}&{\ }\\ \\}\\ \hline
\end{tabular}
\end{center}
\caption{\label{table:sjd} A Haiman table: standard tableaux, their jeu de taquin and
dual equivalence classes}
\end{table}

\subsection{Reading word order?} Further considering
$\Lambda_{{\mathbb O}{\mathbb P}^2}$, we explain
our difficulties in finding a reading
word order for general cominuscule type. In
\cite{Haiman:DE} the reading word of a shape is defined by reading its
entries from left to right and from bottom to top, one row at a time.    
For shapes (respectively,
shifted shapes), \cite{Haiman:DE} gives a short list of 
pairs of reading words such that if $T$ and $U$ are tableaux of size $m=3$ 
(respectively, $m=4$) then $T\equiv_D U$ if and only if the reading words of
$T$ and $U$ appear on this list.  
However, for ${\mathbb {OP}^2}$, we have the following four tableaux:
$$\tableau {  
   &    &   & {} & {} & {} & {} & {}\\
   &    &   & {} & {} & {}\\
   &    & 4 & 5 & {} \\
{} & {} & 1 & 2 &3}
\tableau {  
   &    &   & {} & {} & {} & {} & {}\\
   &    &   & {} & {} & {}\\
   &    & 3& 5 & {} \\
{} & {} & 1 & 2 &4}, \tableau {
   &   &   & 4 & 5 & {} & {} & {}\\
  &   &    & 1 & 2 & 3 \\
  &   &  {}&{} & {} \\
{} & {} & {} & {} & {} }
\tableau {
   &   &   & 3& 5 & {} & {} & {}\\
  &   &    & 1 & 2 & 4\\
  &   &  {}&{} & {} \\
{} & {} & {} & {} & {} }.$$
The first two are dual equivalent while the second two are not.
These pairs of tableaux, however, clearly have the same pairs of
reading words, with respect to the obvious extension of 
the definition in \cite{Haiman:DE} 
or, indeed, with respect to any reading word order defined exclusively
by planar geometry, since the corresponding entries are in the same
relative positions in the two examples.  

The question of a general cominuscule 
description of a reading word order 
is part of the broader question of finding 
a ``semistandard'' theory, together with a ``lattice word'' 
Schubert calculus rule; see, e.g.,~\cite{Stanley, Stembridge} and
the references therein.

\section{Sch\"{u}tzenberger's evacuation involution}\label{evac}

In this section we show how the cominuscule growth diagram approach leads
to a simple proof that M.~P.~Sch\"utzenberger's evacuation is an involution in the
cominuscule setting.  Again, our proofs parallel 
those in S.~Fomin's Appendix~1 to \cite[Chapter~7]{Stanley}.  

The classical evacuation involution appears prominently in
combinatorial representation theory and algebraic geometry; 
see, e.g., \cite{Stembridge2}, and the references therein.
For $T\in {\rm SYT}_{G/P}(\lambda)$, let 
${\widetilde T}$ be obtained by erasing the entry $1$ of $T$
in $\beta(P)$ (the minimal element of $\Lambda_{G/P}$)
and subtracting $1$ from the 
remaining entries. Let
$\Delta(T)={\tt jdt}_{\beta(P)}({\widetilde T})$.
The {\bf evacuation} ${\tt evac}(T)\in {\rm SYT}_{G/P}(\lambda)$ 
is defined by the shape chain
\[\emptyset={\tt shape}(\Delta^{|\lambda|}(T))-
{\tt shape}(\Delta^{|\lambda|-1}(T))-\ldots - {\tt shape}(\Delta^{1}(T)) - {\tt shape}(T).\]

\begin{Theorem}
\label{thm:evac}
${\tt evac}:{\rm SYT}_{G/P}(\lambda)\to {\rm SYT}_{G/P}(\lambda)$
is an involution, i.e., ${\tt evac}({\tt evac}(T))=T$.
\end{Theorem}

For example, if $T=\tableau{&&&{7}\\&&{4}&{6}&{9}\\{1}&{2}&{3}&{5}&{8}} 
\in {\rm SYT}_{G/P}((1,1,2,3,2))$, iterating $\Delta$ gives
$\tableau{&&&{\ }\\&&{5}&{6}&{8}\\{1}&{2}&{3}&{4}&{7}},$
$\tableau{&&&{\ }\\&&{4}&{7}&{\ }\\{1}&{2}&{3}&{5}&{6}},
\tableau{&&&{\ }\\&&{6}&{\ }&{\ }\\{1}&{2}&{3}&{4}&{5}},
\tableau{&&&{\ }\\&&{5}&{\ }&{\ }\\{1}&{2}&{3}&{4}&{\ }},
\tableau{&&&{\ }\\&&{4}&{\ }&{\ }\\{1}&{2}&{3}&{\ }&{\ }},$
$\tableau{&&&{\ }\\&&{\ }&{\ }&{\ }\\{1}&{2}&{3}&{\ }&{\ }},
\tableau{&&&{\ }\\&&{\ }&{\ }&{\ }\\{1}&{2}&{\ }&{\ }&{\ }},
\tableau{&&&{\ }\\&&{\ }&{\ }&{\ }\\{1}&{\ }&{\ }&{\ }&{\ }}$
and hence ${\tt evac}(T)=
\tableau{&&&{9 }\\&&{4}&{7}&{8}\\{1}&{2}&{3}&{5}&{6}}.$ The reader
can check that ${\tt evac}({\tt evac}(T))=T$.

\noindent
\emph{Proof of Theorem~\ref{thm:evac}:} Express each of the tableaux
\[T,\Delta^{1}(T),\ldots, \Delta^{|\lambda|-1}(T),\Delta^{|\lambda|}(T)=\emptyset\] 
as a shape chain and place them right justified in a triangular growth
diagram. In the example above, we have Table~\ref{table:triangle}.
Noting that each ``minor'' of the table whose southwest corner contains a
``$\emptyset$'' is in fact a growth diagram, it follows that the triangular
growth diagram can be reconstructed using the top row and the growth conditions 
of Theorem \ref{thm:gd}. 
Observe that the right column encodes
${\tt evac}(T)$. By the symmetry of growth diagrams,
it follows that applying the above procedure to ${\tt evac}(T)$ 
would give the same triangular growth diagram, after a reflection
across the antidiagonal. Thus the result follows. \qed

\begin{table}[h]
\begin{center}
\begin{tabular}{|l|l|l|l|l|l|l|l|l|l|}
\hline
$\!\emptyset\!$ & $\!(1)\!$ & $\!(1^2)\!$ & $\!(1^3)\!$ & $\!(1,1,2)\!$ & $\!(1,1,2,1)\!$ &
$(1,1,2,2)$ & $(1,1,2,3)$ & $(1,1,2,3,1)$ & $(1,1,2,3,2)$\\ \hline
& $\emptyset$ & $(1)$ & $(1^2)$ & $(1^3)$ & $(1^4)$ & $(1,1,2,1)$ & $(1,1,2,2)$
& $(1,1,2,2,1)$ & $(1,1,2,2,2)$ \\ \hline
& & $\emptyset$ & $(1)$ & $(1^2)$ & $(1^3)$ & $(1,1,2)$ & $(1,1,2,1)$ &
$(1,1,2,1,1)$ & $(1,1,2,2,1)$ \\ \hline
& & & $\emptyset$ & $(1)$ & $(1^2)$ & $(1^3)$ & $(1^4)$ & $(1^5)$ &
$(1,1,2,1,1)$ \\ \hline
& & & & $\emptyset$ & $(1)$ & $(1^2)$ & $(1^3)$ & $(1^4)$ & $(1,1,2,1)$ \\ \hline
& & & & & $\emptyset$ & $(1)$ & $(1^2)$ & $(1^3)$ & $(1,1,2)$ \\ \hline
& & & & & & $\emptyset$ & $(1)$ & $(1^2)$ & $(1^3)$ \\ \hline
& & & & & & & $\emptyset$ & $(1)$ & $(1^2)$ \\ \hline
& & & & & & & & $\emptyset$ & $(1)$ \\ \hline
& & & & & & & & & $\emptyset$ \\ \hline
\end{tabular}
\end{center}
\caption{\label{table:triangle} A triangular growth diagram, for the proof of Theorem~\ref{thm:evac}.}
\end{table}

\section{Cartons}

The goal of this section is to extend the main result of
\cite{Thomas.Yong:S3} to the cominuscule setting. Our description
of the rule closely parallels the one for the original rule from 
our earlier paper.

\subsection{Statement of the rule} 

Let $\lambda$, $\mu$, and $\nu$ be shapes in $\mathbb Y_{G/P}$, such that
$|\lambda|+|\mu|+|\nu|=|\Lambda_{G/P}|$.  (When this condition is
not satisfied, $c_{\lambda,\mu,\nu}(G/P)$ is necessarily 0.)
 
\begin{figure}[h]
\unitlength=.32mm
\begin{picture}(120,158)
\thicklines
\put(60,10){\line(1,1){50}}
\put(60,25){\line(1,1){50}}
\put(60,41){\line(1,1){50}}
\put(60,10){\line(-1,1){50}}
\put(60,10){\line(0,1){50}}
\put(110,60){\line(0,1){50}}

\put(10,110){\line(1,1){50}}
\put(25,95){\line(1,1){50}}
\put(42,78){\line(1,1){50}}
\put(70,20){\line(0,1){50}}
\put(80,30){\line(0,1){50}}
\put(90,40){\line(0,1){50}}
\put(100,50){\line(0,1){50}}
\put(10,60){\line(0,1){50}}
\put(25,45){\line(0,1){50}}
\put(42,28){\line(0,1){50}} 

\put(110,110){\line(-1,1){50}}
\put(100,100){\line(-1,1){50}}

\put(60,60){\line(-1,1){50}}
\put(60,41){\line(-1,1){50}}
\put(60,25){\line(-1,1){50}}
\put(90,90){\line(-1,1){50}}
\put(80,80){\line(-1,1){50}}
\put(70,70){\line(-1,1){50}}

\put(60,60){\line(1,1){50}}
\put(57,-1){$\emptyset$}
\put(113,55){$T_\mu$}
\put(54,68){$T_\lambda$}
\put(-2,55){$T_\nu$}
\put(57,162){$\Lambda_{G/P}$}
\end{picture}
\caption{\label{fig:growth} Theorem~\ref{thm:main2} 
calculates $c_{\lambda,\mu,\nu}(G/P)$ by
assigning Young diagrams to the vertices of the six faces}
\end{figure}
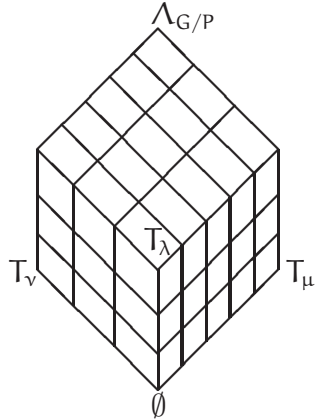

Figure~\ref{fig:growth} depicts a {\bf carton}.  This is a 
$|\lambda|\times|\mu|\times |\nu|$ box, with a 
grid whose squares are $1\times 1$ drawn on each of the six faces.  
One vertex of the box is labelled $\emptyset$, and its opposite vertex
is labelled $\Lambda_{G/P}$.  
A {\bf carton filling} assigns a Young diagram
to each vertex of the grid so that shapes increase one box at a time while
moving away from the $\emptyset$, so that for any $2\times 2$ subgrid
$\begin{matrix}
\alpha & - & \beta\\
| & & | \\
\gamma & - &\delta
\end{matrix}$ the {\bf Fomin growth conditions} of Theorem~\ref{thm:gd} hold.
Note that there are vertices of the grid which lie on edges of the box, and
thus participate in $2\times 2$ subgrids on more than one face.

Fix a choice of standard tableaux $T_{\lambda},T_{\mu}$ and
$T_{\nu}$ of respective shapes $\lambda,\mu$ and $\nu$. 
Initialize the edges $\emptyset-T_{\lambda}$, $\emptyset-T_{\mu}$
and $\emptyset-T_{\nu}$ 
with the shape chains for the corresponding tableaux.
Let ${\tt CARTONS}_{\lambda,\mu,\nu}(G/P)$ be
all carton fillings with the above 
initial data.

\begin{Theorem}
\label{thm:main2}
For cominuscule $G/P$, 
\[c_{\lambda,\mu,\nu}(G/P)=2^{{\tt shortroots}(\Lambda_{G/P})-{\tt shortroots}(\nu)-{\tt shortroots}(\lambda)-{\tt shortroots}(\mu)}\#{\tt CARTONS}_{\lambda,\mu,\nu}(G/P).\] 
\end{Theorem}

\begin{figure}[h]
\unitlength=.33mm
\begin{picture}(200,360)
\thicklines
\put(80,30){\line(1,1){150}}
\put(80,30){\line(-1,1){50}}
\put(80,30){\line(0,1){100}}
\put(80,130){\line(1,1){150}}
\put(230,180){\line(0,1){100}}
\put(30,80){\line(0,1){100}}
\put(30,180){\line(1,1){150}}
\put(30,180){\line(1,-1){50}}
\put(180,330){\line(1,-1){50}}
\put(80,80){\line(1,1){150}}
\put(80,80){\line(-1,1){50}}
\put(130,80){\line(0,1){100}}
\put(130,180){\line(-1,1){50}}
\put(180,130){\line(0,1){100}}
\put(180,230){\line(-1,1){50}}
\put(77,10){$\emptyset$}
\put(130,60){$\tableau{{\ }}$}
\put(180,110){$\tableau{{\ }&{\ }}$}
\put(235,170){$\tableau{&{3}\\{1 } &{2 }}=\mu$}
\put(235,230){$\tableau{&{\ }\\{\ }&{ \ }&{\  }}$}
\put(235,280){$\tableau{&{\ }&{\ }\\{\ }&{\ }&{\ }}$}
\put(85,70){$\tableau{{\ }}$}
\put(135,120){$\tableau{{\ }&{\ }}$}
\put(185,170){$\tableau{{\ }&{\ }&{\  }}$}
\put(85,120){$\tableau{{1}&{2}}$}
\put(85,105){$=\lambda$}
\put(135,175){$\tableau{&{\ }\\{\ }&{ \ }}$}
\put(185,230){$\tableau{&{\ }\\{\ }&{ \ }&{ \  }}$}
\put(-12,70){$\nu=\tableau{{1 }}$}
\put(0,120){$\tableau{{ \ }&{\ }}$}
\put(-12,170){$\tableau{{\ }&{\ }&{\ }}$}
\put(45,250){$\tableau{&{\ }\\{\ }&{\  }&{\ }}$}
\put(95,300){$\tableau{&{\ }&{\ }\\{\ }&{\ }&{\  }}$}
\put(162,360){$\tableau{&&{ \ }\\&{\ }&{ \ }\\{\ }&{\ }&{ \ }}=\Lambda_{G/P}$}
\end{picture}
\caption{\label{fig:easyex} The ``front'' three faces of the (unique)
carton filling for $c_{(1,1,0),(1,2,0),(1,0,0)}(LG(3,6))=2^{3-1-1}\cdot 1$. 
The tableaux $T_{\lambda}, T_{\mu}$ and
$T_{\nu}$ are as shown.}
\end{figure}
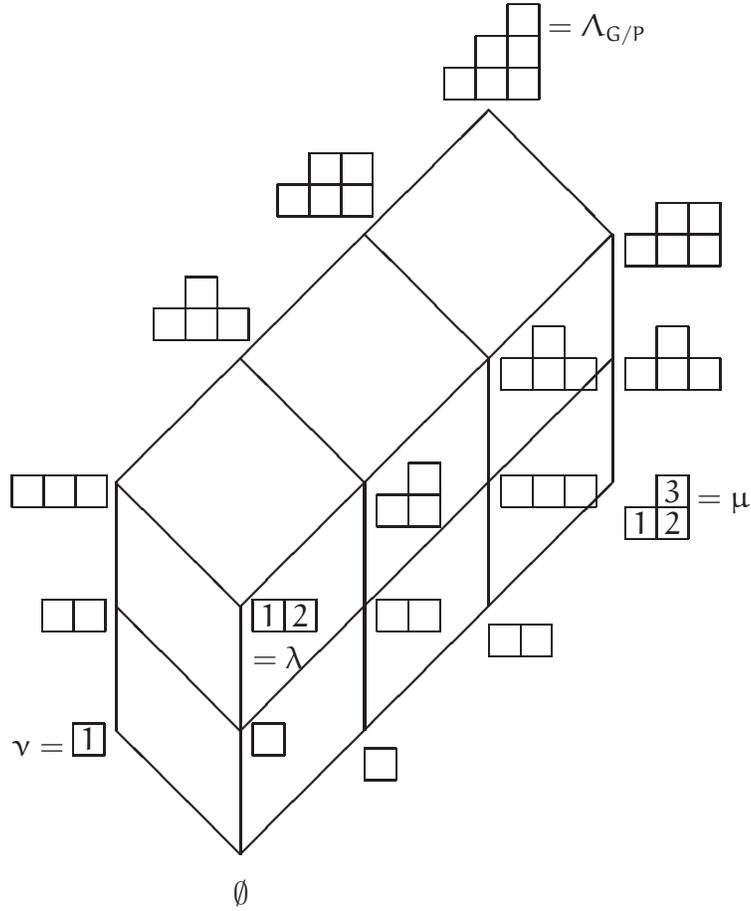

This rule manifests bijections between 
${\tt CARTONS}_{\lambda,\mu,\nu}(G/P)$
and ${\tt CARTONS}_{\alpha,\beta,\gamma}(G/P)$ for any permutation
$(\alpha,\beta,\gamma)$ of $(\lambda,\mu,\nu)$. 
In Figure~\ref{fig:easyex} we give an example of Theorem~\ref{thm:main2}.

\subsection{The proof}
The proof in the Grassmannian case is given in \cite{Thomas.Yong:S3}. 
It carries over to the cominuscule setting  
using the tools developed for that setting in 
the previous sections (namely: dual equivalence, growth diagrams, and
evacuation).  Since the proof in the cominuscule case is the same as
in the Grassmannian case, we do not give all the details, as the interested
reader will have no trouble filling them in from \cite{Thomas.Yong:S3}.


Let $\alpha \in \mathbb Y_{G/P}$.  
Recall that we write $\alpha^\vee$ for the shape obtained by taking any
standard tableau of shape $\alpha$ and reverse rectifying it as far as 
possible.  It is easy to see that if 
$\alpha \supset \beta$, then $\beta^\vee \supset 
\alpha^\vee$.  Thus $\alpha \rightarrow \alpha^\vee$ is an 
anti-automorphism of $Y_{G/P}$.  It therefore induces an isomorphism from the
join-irreducibles of $Y_{G/P}$ to the join-irredicibles of the poset of
reverse shapes.  This in turn induces an anti-automorphism of 
$\Lambda_{G/P}$.  We call this anti-automorphism {\tt rotate}, because in
the Grassmannian case it amounts to rotation by 180 degrees.  An explicit
description of ${\tt rotate}$ in each of the cominuscule cases can be found
in \cite[Section~2.2]{Thomas.Yong}, whose equivalence with the definition 
we have given here is 
\cite[Proposition~4.6]{Thomas.Yong}.
(We alluded earlier to the easily checked
fact that $\Lambda_{G/P}$ is self-dual, and thus admits
\emph{some} antiautomorphism.  However, sometimes there is more than
one anti-automorphism, and in those cases, it is important to use the correct
one, defined as above.)

Given $T$ a tableau of shape $\alpha$, we denote by
${\tt rotate(T)}$ the tableau of shape $\Lambda_{G/P}/\alpha^\vee$ obtained 
by putting the label from box $x$ into the box ${\tt rotate}(x)$.  

Given $T\!\in\! {\rm SYT}(\alpha)$ for a straight shape $\alpha$, define
${\widetilde T}\!\in\! {\rm SYT}({\tt rotate}(\alpha))$ by computing 
${\tt evac}(T)\in {\rm SYT}(\alpha)$, replacing entry $i$ with $|\alpha|-i+1$
throughout and applying ${\tt rotate}$.

The following fact extends \cite[Lemma~2.1]{Thomas.Yong:S3} with the same proof, given our 
definition of {\tt rotate} above and Corollary \ref{cor:proctor}.  

\begin{Lemma}
Let $\alpha,\beta,\gamma\in {\mathbb Y}_{G/P}$ and let
$T_{\beta}\in {\rm SYT}(\beta), T_{\gamma^{\vee}/\alpha}\in 
{\rm SYT}(\gamma^{\vee}/\alpha)$
be tableaux satisfying 
${\tt rectification}(T_{\gamma^{\vee}/\alpha})=T_{\beta}$.
Then 
\[{\tt revrectification}(T_{\gamma^{\vee}/\alpha})={\tt revrectification}(T_{\beta})={\widetilde T}_{\beta}.\]
\end{Lemma}

As in \cite[Corollary~2.2]{Thomas.Yong:S3}, we have:

\begin{Corollary}
\label{cor:moreinit}
Fix a carton filling. The face joining the edges assigned 
the shape chains for $T_\lambda$ and $T_\mu$, necessarily has assigned to
its uninitialized corner the shape $\nu^\vee$.  Similarly, the face
joining the edges $T_\lambda$ and $T_\nu$, has assigned to its uninitialized
corner the shape $\mu^\vee$, and the face joining the edges
$T_\mu$ and $T_\nu$, has assigned to its uninitialized corner (the
corner not visible in Figure~\ref{fig:growth}) the shape $\lambda^\vee$.  
Thus, we can refer to the
edges $\lambda^{\vee}-\Lambda$, $\mu^{\vee}-\Lambda$ and $\nu^{\vee}-\Lambda$.
These edges are necessarily assigned 
the shape chains of ${\widetilde T}_{\lambda}, {\widetilde T}_{\mu}$
and ${\widetilde T}_{\nu}$ respectively.
\end{Corollary}

Thus by Corollary~\ref{cor:moreinit}, it makes sense to refer to a
face by its corner vertices. 
Note any carton filling
gives a growth diagram on the face 
$\emptyset-\mu-\nu^{\vee}-\lambda$ for which the edge $\lambda-\nu^{\vee}$
is a standard tableau of shape $\nu^{\vee}/\lambda$ rectifying
to $T_{\mu}$. By Theorem~\ref{thm:originalmain}, 
fillings of this face count 
$2^{{\tt shortroots}(\mu)-{\tt shortroots}(\nu^{\vee}/\lambda)}c_{\lambda,\mu,\nu}(G/P)$. 

Conversely, if we start with a filling of the 
$\emptyset-\mu-\nu^\vee-\lambda$ face, then it is straightforward to
use the Fomin growth
conditions and Corollary~\ref{cor:moreinit} to show that there is at most
one way to extend this filling to a filling of the entire carton.  The proof
that there is exactly one way to extend the filling  
follows exactly as in
\cite[Section~2.2]{Thomas.Yong:S3}, where references to the 
growth diagram encoding
of evacuation from S.~Fomin's \cite[Appendix~1]{Stanley} are replaced by
the cominuscule generalization given in Section \ref{evac} above.  
Finally, one observes that
$2^{{\tt shortroots}(\nu^{\vee}/\lambda)-{\tt shortroots}(\mu)}=
2^{{\tt shortroots}(\Lambda_{G/P})-{\tt shortroots}(\nu)-{\tt shortroots}(\lambda)-{\tt shortroots}(\mu)}$.\qed

\section*{Acknowledgements}
We thank Sami Assaf and Mark Haiman 
for their informative explanations of dual equivalence, 
during an NSF RTG supported visit of AY to UC Berkeley in October 2006. In particular,
Assaf's work on dual equivalence and Macdonald polynomials
made us curious about how dual equivalence might apply to 
Schubert calculus. We would also like to thank Stephen Griffeth, 
Oliver Pechenik, Robert Proctor, Kevin Purbhoo, Victor Reiner, Muge Taskin and 
two anonymous referees for helpful communications.

HT was supported by an NSERC Discovery Grant, and would like to thank the Fields Institute for its hospitality. AY was
partially supported by NSF grants and an NSERC Postdoctoral
Fellowship held at the Fields Institute, Toronto. This text was 
also completed while AY was a Helen Corley Petit Scholar at UIUC.

Finally, we would like to thank the organizers of the International
Seasonal Institute on Schubert calculus (MSJ-SI 2012) in Osaka, Japan for providing the
encouraging environment and circumstances to complete this
text in its present form.


\begin{thebibliography}{ClThYo12}
\bibitem[BeSoSt96]{Benkart} G.~Benkart, F.~Sottile and J.~Stroomer,
\emph{Tableau switching: Algorithms and applications,} J.~Combin.~Theory
Ser.~A {\bf 76}(1)(1996), 11--43.

\bibitem[BuRa12]{Buch.Ravikumar} A.~Buch and V.~Ravikumar,
\emph{Pieri rules for the $K$-theory of cominuscule 
Grassmannians}, J.~Reine Angew. Math. {\bf 668}(2012),
109--132.
\bibitem[BuSa13]{Buch.Samuels} A.~Buch and M.~Samuels,
\emph{$K$-theory of minuscule varieties}, preprint, 2013. 
\textsf{arXiv:1306.5419}
\bibitem[ChPe12]{Chaput.Perrin} P.~E.~Chaput and N.~Perrin, 
\emph{Towards a Littlewood-Richardson rule for Kac-Moody homogeneous spaces}, J.~Lie~Theory {\bf 22}(2012), 17--80.
\bibitem[ClThYo12]{Clifford.Thomas.Yong} E.~Clifford, H.~Thomas
and A.~Yong, \emph{$K$-theoretic Schubert calculus for 
$OG(n,2n+1)$ and jeu de taquin for shifted increasing tableaux},
J.~Reine Angew. Math., to appear, 2012. 
\bibitem[Ha92]{Haiman:DE} M.~Haiman, \emph{Dual equivalence with applications,
including a conjecture of Proctor}, Discrete Math.~{\bf 99}(1992), 79--113.
\bibitem[Pe12]{Pechenik} O.~Pechenik, \emph{Cyclic sieving of
increasing tableaux and small Schr\"{o}der paths}, preprint, 2012.
\textsf{arXiv:1209.1355}
\bibitem[Pr84]{Proctor:euro} R.~Proctor, \emph{Bruhat lattices, plane
partition generating functions, and minuscule representations}, Europ.~J.~Combinatorics, {\bf 5}(1984), 331--350.
\bibitem[Pr04]{Proctor} \bysame, \emph{d-Complete posets generalize Young diagrams for the jeu de taquin property}, preprint, 2004, available at
\textsf{http://www.math.unc.edu/Faculty/rap/}

\bibitem[Sa87]{sagan} B.~E.~Sagan, \emph{Shifted tableaux, Schur~$Q-$ functions, and a conjecture
of Stanley}, J.~Combin.~Theory
Ser.~A {\bf 45}(1987), 62--103.
\bibitem[Sc77]{Schutzenberger} M.-P.~Sch\"{u}tzenberger, \emph{Combinatoire et repr\'{e}sentation du 
groupe sym\'{e}trique} (Actes Table Ronde CNRS, Univ.~Louis-Pasteur Strasbourg, Strasbourg, 1976),  
pp. 59--113. Lecture Notes in Math., Vol. 579, Springer, Berlin, 1977. 

\bibitem[St99]{Stanley} R.~P.~Stanley, \emph{Enumerative Combinatorics, Volume~2} (with an appendix by S.~Fomin), Cambridge University Press, 1999.

\bibitem[St09]{stanley:recent} \bysame, \emph{Promotion and evacuation},
Electronic~J.~Combin. {\bf 16}(2)(2009), R9.  


\bibitem[St89]{Stembridge} J.~Stembridge, \emph{Shifted tableaux and the 
projective representations of symmetric groups}, Adv.~Math.~{\bf 74}(1989),
87--134.

\bibitem[St96]{Stembridge2} \bysame,
\emph{Canonical Bases and Self-Evacuating Tableaux},
    Duke Math. J. {\bf 82}(1996), 585--606.


\bibitem[ThYo08]{Thomas.Yong:S3} H.~Thomas and A.~Yong, \emph{An $S_3$-symmetric Littlewood-Richardson rule},
Math.~Res.~Lett. {\bf 15}(2008), no. 5, 1027--1037.

\bibitem[ThYo09a]{Thomas.Yong} \bysame, \emph{A combinatorial
rule for (co)minuscule Schubert calculus}, Adv.~Math.~{\bf 222}(2009), no. 2, 596--620.

\bibitem[ThYo09b]{Thomas.Yong:K} \bysame, \emph{A jeu de taquin theory
for increasing tableaux, with applications to $K$-theoretic Schubert calculus}, Algebra Number Theory {\bf 3}(2009), no. 2, 121--148.

\bibitem[Wo84]{worley} D.~Worley, \emph{A theory of shifted Young tableau}, Ph.~D.~thesis, M.~I.~T., 1984, available at {\tt http://hdl.handle.net/1721.1/15599} .
\end{thebibliography}
\end{document}